\documentclass[11pt]{article}
\usepackage{amsmath, amsfonts, amssymb, amsthm,epsfig}

\let\BFseries\bfseries\def\bfseries{\BFseries\mathversion{bold}} 
\DeclareMathSymbol{\leqslant}{\mathalpha}{AMSa}{"36} 
\DeclareMathSymbol{\geslant}{\mathalpha}{AMSa}{"3E} 
\DeclareMathSymbol{\eset}{\mathalpha}{AMSb}{"3F}     

\newcommand{\ind}{1\hspace{-0.098cm}\mathrm{l}}

\newcommand{\eps}{\varepsilon}

\newcommand{\IR}{\mathbb{R}}

\newcommand{\be}{\begin{eqnarray*}}
\newcommand{\ee}{\end{eqnarray*}}
\newcommand{\ben}{\begin{eqnarray}}
\newcommand{\een}{\end{eqnarray}}

\theoremstyle{plain}
\newtheorem{thm}{Theorem}
\newtheorem{lem}[thm]{Lemma}

\newtheorem{cor}[thm]{Corollary}

\theoremstyle{definition}

\newtheorem{remark}[thm]{Remark}

\renewenvironment{proof}[1][] {\smallskip \noindent {\bf Proof#1.} }{\hspace*{\fill}$\square$\medskip\par}
\def\P{{\bf {\mathbb{P}}}}
\newcommand{\pr}[1]{\P\left[#1\right]}

\def\E{\mathbb{E}}

\def\R{\IR}
\def\d{{\rm d}}

\begin{document}
\title{On the one-sided exit problem for fractional Brownian motion}
\author{Frank Aurzada}
\date{\today}
\maketitle
\begin{abstract}
 We consider the one-sided exit problem for fractional Brownian motion (FBM), which is equivalent to the question of the distribution of the lower tail of the maximum of FBM on the unit interval. We improve the bounds given by Molchan (1999) and shed some light on the relation to the quantity $I$ studied there.
\end{abstract}



\section{Introduction}
This paper is concerned with the so-called one-sided exit problem for stoch\-astic processes: If $(X(t))_{t\geq 0}$ is a real-valued stochastic process, we want to find the asymptotic rate of the following probability:
\begin{equation} \label{eqn:rate}
F(T):=\pr{ \sup_{0\leq t\leq T} X(t) \leq 1},\qquad\text{when $T\to\infty$.}
\end{equation}
This problem arises in a number of contexts, the most important of which is the relation to Burgers equation with random initial data (see e.g.\ \cite{bertoin,frisch}). Further applications concern pursuit problems, relations to random polynomials, and polymer models. We refer to \cite{lishao2004aop} for more background information and links to further literature.

In the context of Gaussian process, there seem to be very few results concerning the asymptotic rate of (\ref{eqn:rate}). The precise rate is known only for Brownian motion and integrated Brownian motion and very few very particular Gaussian processes.

The aim of this paper is to study the rate in (\ref{eqn:rate}) for fractional Brownian motion. Fractional Brownian motion (FBM) $X$ is a centered Gaussian process with covariance
$$
\E X(t) X(s) = \frac{1}{2} \left( |t|^{2H} + |s|^{2H} - |t-s|^{2H}\right), \qquad t,s\in\R,
$$
where $0<H<1$ is the so-called Hurst parameter. It is well-known that FBM is self-similar with index $H$ and has stationary increments.

The question of the exit probability (\ref{eqn:rate}) for FBM has been investigated by Sina\u\i\, \cite{sinai1997} and Molchan (\cite{molchan1999,molchan1999b,molchan2008}). The most precise result concerning the asymptotics in (\ref{eqn:rate}) for fractional Brownian motion, \cite{molchan1999}, states that
\begin{equation}\label{eqn:molchan}
 T^{-(1-H)}\, e^{-k \sqrt{\log T}} \leq \pr{ \sup_{0\leq t\leq T} X(t) \leq 1} \leq T^{-(1-H)}\, e^{+k \sqrt{\log T}},
\end{equation}
for some positive constant $k$ and $T$ large enough.

In the physics literature, this result is used in the sense of $\approx T^{-(1-H)}$, disregarding the loss factors $e^{\pm k \sqrt{\log T}}$. We stress that already proving (\ref{eqn:molchan}) is highly non-trivial and that presently there is no approach to obtain the precise order of this probability. We mentioned that, beyond the classical results for such particular processes as Brownian motion or integrated Brownian motion, there is no theory to obtain even the polynomial term. Due to this lack of theory, even for simple-looking estimates rather involved calculations are needed, see e.g.\ (\ref{eqn:difficultsiesstimate}) below.

In this paper, we give the following improvement of (\ref{eqn:molchan}).

\begin{thm} \label{thm:main}
There is a constant $c>0$ such that, for large enough $T$, we have
$$
T^{-(1-H)}\, (\log T)^{-c} \leq \pr{ \sup_{0\leq t\leq T} X(t) \leq 1} \leq T^{-(1-H)}\, (\log T)^{+c} .
$$
\end{thm}

Before giving the proofs of the lower and upper bound in Sections~\ref{sec:lower} and~\ref{sec:upper}, respectively, we will give some comments.

Molchan \cite{molchan1999} related the problem of finding the asymptotics of $F$ to the quantity
\begin{equation} \label{eqn:defnI}
I(T):=\E\left[  \left( \int_0^T e^{X(u)} \d u \right)^{-1} \right];
\end{equation}
and he is even able to determine the strong asymptotic rate of $I$. However, when passing over from $I$ to $F$, the slowly varying terms $e^{\pm k \sqrt{\log T}}$ appear. This is essentially due to a change of measure argument: if $g$ is a function in the reproducing kernel Hilbert space of $X$ then the asymptotic rates for the exit problems of $X+g$ and $X$, respectively, differ at most by $e^{\pm k \sqrt{\log T}}$, cf.\ \cite{aurzadadereich2009}, Proposition~3.1.

A main goal of this paper is to shed light on the relation between $I$ and $F$. Heuristically, it is clear that those paths of $X$ that remain below $1$ until $T$ will escape to $-\infty$ rather rapidly and thus give a major contribution to $I$. Vice versa, those paths that do not remain below $1$ until $T$ will tend to be near or above zero for a positive fraction of time and thus do not give much contribution to $I$.

Our proofs will make an effort to understand this relation -- beyond a heuristic level. The proof of the upper bound in Theorem~\ref{thm:main} is based on seeing $I(T)$ as an exponential integral. The proof of the lower bound in Theorem~\ref{thm:main} selects some paths in the expectation in (\ref{eqn:defnI}) that give a relevant contribution.

The constant $c$ appearing in the theorem can be specified. For the lower bound one can choose any $c>1/(2H)$, for the upper bound any $c>2/H-1$. However, we do not conjecture optimality of either of the constants. In fact, our proofs make it plausible that even $F(T)\approx T^{-(1-H)}$.

Due to the self-similarity of fractional Brownian motion, our main result immediately translates into a result for the lower tail of the maximum of fractional Brownian motion.

\begin{cor}
There is a constant $c>0$ such that, for small enough $\eps$, we have
$$
\eps^{(1-H)/H} \, |\log \eps|^{-c} \leq \pr{ X^*_1 \leq \eps} \leq \eps^{(1-H)/H}\, |\log \eps|^{+c},
$$
where $X^*_1:=\sup_{0\leq t\leq 1} X(t)$.
\end{cor}

\section{Lower bound} \label{sec:lower}
Before proving the lower bound in Theorem~\ref{thm:main}, we explain the main line of thought. It shows that the quantity $I$ is indeed a natural object in the study of one-sided exit probabilities (\ref{eqn:rate}), even beyond FBM.

In the following, we use $X^*_1:=\sup_{u\in[0,1]} X(u)$ and $u^*$ a point in $[0,1]$ where $X(u^*)=X^*_1$.

The self-similarity of $X$ implies that
\begin{eqnarray*}
I(T)&=&\E\left[  \left( \int_0^1 e^{X(T u)} T \d u \right)^{-1} \right]=\E\left[  \left( \int_0^1 e^{T^H X(u)} T \d u \right)^{-1} \right]\\
&=&\E\left[  \left( \int_0^1 e^{-T^H (X(u^*)-X(u))} T \d u  \right)^{-1} e^{-T^H X^*_1} \right].
\end{eqnarray*}

The path of the process $X$ is H\"older continuous with H\"older exponent $\gamma<H$. In particular, it is not H\"older continuous with exponent $H$. However, suppose for a moment that $|X(t)-X(s)|\sim S |t-s|^H$ for $|t-s|\to 0$. Then the above term behaves asymptotically as
\begin{eqnarray*}
I(T)&\approx&\E\left[  \left( \int_{u^*-\eps}^{u^*+\eps} e^{-T^H S |u^* -u|^H} T \d u  \right)^{-1} e^{-T^H X^*_1} \right] \\
&\approx&\E\left[  \left( \int_{-\infty}^{\infty} e^{-S |x|^H} \d x  \right)^{-1} e^{-T^H X^*_1} \right] \\
&=&  \left( \int_{-\infty}^{\infty} e^{-S |x|^H} \d x  \right)^{-1} \E\left[ e^{-T^H X^*_1} \right] \\
&=&  c\, \E\left[ e^{-T^H X^*_1} \right].
\end{eqnarray*}
Now, via Tauberian theorems, the behaviour of $\E\left[ e^{-T^H X^*_1} \right]$ as $T\to\infty$ is related to the one-sided exit problem $\pr{ X^*_1 \leq \eps}$ as $\eps\to 0$. However, by the self-similarity of $X$, we have
$$
\pr{ \sup_{0\leq s \leq 1} X(s) \leq \eps} \sim c' \eps^{\theta/H}
\quad\Leftrightarrow\quad
\pr{ \sup_{0\leq s \leq T} X(s) \leq 1} \sim c  T^{-\theta},
$$
which brings us back to our original problem.

Of course, fractional Brownian motion does not satisfy $|X(t)-X(s)|\sim S |t-s|^H$, so that the above calculations are just heuristics. However, the idea can be turned into a formally correct proof of the lower bound for $F$.

\begin{proof}[ of the lower bound in Theorem~\ref{thm:main}]
{\it Step 1:} The crucial inequality.

Let $H/2<\gamma<H$. Fix $a$ such that $a>2/H>1/\gamma$ and $\gamma<H-1/a$. Due to the stationarity of increments, it is clear that fractional Brownian motion satisfies
$$
\left( \E |X(t)-X(s)|^a \right)^{1/a} = C(a) |t-s|^H,\qquad t,s\geq 0.
$$
A close analysis of the constant $C(a)$ shows that $C(a) \leq c a^{1/2}$, as $a\to \infty$.
Therefore, we have
\begin{equation} \label{eqn:momentcon}
\E |X(t)-X(s)|^a \leq c^a a^{s a} |t-s|^{1+aH-1},\qquad t,s\in[0,1],
\end{equation}
where $s:=1/2$ (for readability), $c>0$, and $aH-1>0$.

By the well-known Kolmogorov theorem, this implies that $X$ has H\"older continuous paths of order $\gamma$. An extension of Kolmogorov's theorem (see \cite{scheutzow}, Lemma~2.1) implies even more: namely, an estimate for the modulus of H\"older continuity. Concretely, from (\ref{eqn:momentcon}) we can infer that, for any $0<\eps\leq 1$,
$$
|X(t)-X(s)| \leq S \eps^\gamma,\qquad\text{for all $t,s\in[0,1]$ with $|t-s|\leq \eps$},
$$
where $S$ is a random variable with
$$
\E S^a \leq \frac{2^a}{(1-2^{-\gamma})^a} \, \frac{(c a^s)^a}{2^{(aH-1)-a \gamma}-1}
  \leq 
 \frac{(d a^s)^{a}}{2^{(aH-1)-a \gamma}-1}.
$$
Let us now mimic the heuristics presented before the proof: For $0<\eps\leq 1$,
\begin{eqnarray*}
 I(T)
 &=&\E\left[  \left( \int_0^1 e^{-T^H (X^*_1 - X(u))} T \d u \right)^{-1} e^{-T^H X^*_1}\right] \\
 &\leq &\E\left[  \left( \int_{(u^*-\eps)\vee 0}^{(u^*+\eps)\wedge 1} e^{-T^H (X(u^*) - X(u))} T \d u \right)^{-1} e^{-T^H X^*_1}\right]\\
 &\leq &\E\left[  \left( \int_{(u^*-\eps)\vee 0}^{(u^*+\eps)\wedge 1} e^{-T^H S \eps^\gamma} T \d u \right)^{-1} e^{-T^H X^*_1}\right]\\
 &= &\E\left[ \eps^{-1} e^{T^H S \eps^\gamma} T^{-1} e^{-T^H X^*_1}\right].
\end{eqnarray*}
Setting $\eps:=(T^H S)^{-1/\gamma}\wedge 1$ yields $T^H S \eps^\gamma\leq 1$ and $\eps^{-1}\leq (T^H S)^{1/\gamma} +1$. Thus,
\begin{align*}
I(T) & \leq \E\left[ (T^{H/\gamma} S^{1/\gamma} + 1) e^1 e^{-T^H X^*_1}\right] T^{-1} \\
 & = e \left( \E\left[ T^{H/\gamma} S^{1/\gamma} e^{-T^H X^*_1}\right] + \E\left[ e^{-T^H X^*_1}\right] \right)T^{-1}.
\end{align*}
For simplicity, we set $g(T):=\E \left[ e^{-T^H X^*_1}\right]$. Now we use H\"older's inequality ($1/p+1/q=1$, $p,q>1$) in the first term to get
$$
\frac{I(T)}{e} \leq \E\left[ S^{p/\gamma} \right]^{1/p} \E \left[ e^{-q T^H X^*_1}\right]^{1/q} T^{-1+H/\gamma} +T^{-1} g(T).
$$
Setting $p/\gamma:=a$ and using the estimate for the $a$-th moment of $S$ (and $q\geq 1$) we obtain:
$$
\frac{I(T)}{e}  \leq \left( \frac{(d a^s)^a}{2^{(aH-1)-a \gamma}-1} \right)^{1/(a\gamma)} \E \left[ e^{-T^H X^*_1}\right]^{1/q} T^{-1+H/\gamma}+ T^{-1} g(T).
$$
We rewrite this as follows:
\begin{equation} \label{eqn:crucial}
k \, I(T)  \leq \frac{a^{s/\gamma}}{(2^{a(H-\gamma)-1}-1)^{1/(a\gamma)}} \,g(T)^{1-1/(a\gamma)} T^{-1+H/\gamma}+ T^{-1} g(T).
\end{equation}
Here $k$ is some constant depending only on $H$; and this inequality holds for all $T>0$, $\gamma<H$, $a>2/H$, and $\gamma<H-1/a$.

{\it Step 2:} We use Molchan's result for $I(T)$ to conclude.

Fix $0<\delta<1$. We choose
\begin{equation} \label{eqn:choiceofaa}
a:=\log T \cdot (\log \log T)^{-\delta},\qquad \gamma:=H-2/a.
\end{equation}
Then, by (\ref{eqn:crucial}), for $T$ large enough,
\begin{align*}
k\, I(T) & \leq \frac{a^{s/\gamma}}{(2^{2-1}-1)^{1/(a\gamma)}}\, g(T)^{1-1/(a\gamma)} T^{-1+1+2/(a\gamma)} + T^{-1} g(T)\\ & =  a^{s/\gamma} g(T)^{1-1/(a\gamma)} T^{2/(a\gamma)} + T^{-1} g(T)\\
& \leq 2 a^{s/\gamma} g(T)^{1-1/(a\gamma)} T^{2/(a\gamma)}.
\end{align*}
Rewriting the inequality we get
$$
\left( k' I(T) a^{-s/\gamma} T^{-2/(a\gamma)} \right)^{1/(1-1/(a\gamma))} \leq g(T).
$$

We know from Molchan (Statement~1 in \cite{molchan1999}) that $I(T)\geq c T^{-(1-H)}$ for $T$ large enough. Then the left-hand side becomes:
\begin{align*}
& \exp\left( (1+O(a^{-1}))( \log  k'' -(1-H) \log T - \frac{s}{\gamma} \log a - \frac{2}{a\gamma} \log T \right)
\\
& = \exp\left( -(1-H) \log T - \frac{s}{\gamma} \log \log T - o(\log \log T) \right)
\\
& =  \exp\left( -(1-H) \log T - \frac{s}{H} \log \log T +o(\log \log T) \right)
\\
&= T^{-(1-H)} (\log T)^{-s/H+o(1)}.
\end{align*}
The remainder of the proof is clear due to Tauberian theorems (see \cite{bgt}, Corollary~8.1.7):
$$
g(T)\geq T^{-(1-H)} (\log T)^{-s/H+o(1)}
$$
implies the asymptotic lower bound for $F(T)$.
\end{proof}

\begin{remark}
We remark that we have used few properties of fractional Brownian motion. In fact, only (\ref{eqn:momentcon}) and the self-similarity is used, not even Gaussianity.
\end{remark}

\section{Upper bound} \label{sec:upper}
The main idea of the proof of the lower bound is to restrict the expectation in (\ref{eqn:defnI}) to a set of paths where the integral can be estimated.

In the proof we have to distinguish the cases of positively ($H\geq 1/2$) and negatively ($H<1/2$) correlated increments. The latter case is more involved but contains the same main idea.

{\bf Proof of the lower bound in Theorem~\ref{thm:main} for the case $H\leq 1/2$.}
Let $\kappa>1$ and define
$$
\phi(t):=\begin{cases}
          1  & 0\leq t \leq (\kappa\log T)^{1/H},\\
	  - \kappa \log T & (\kappa\log T)^{1/H} \leq t \leq T.
         \end{cases}
$$
Clearly,
\begin{align}
I(T) & =\E \left[\left(\int_0^T e^{ X(u)} \, \d u \right)^{-1} \right] \notag \\
 & \geq \E  \left[ \ind_{\{\forall 0\leq t\leq T : X(t) \leq \phi(t)\}} \cdot \left(\int_0^T e^{ X(u)} \, \d u \right)^{-1} \right] \notag \\
 & \geq \pr{\forall 0\leq t\leq T : X(t) \leq \phi(t)}\cdot \left(\int_0^T e^{ \phi(u)} \, \d u \right)^{-1} \notag \\
 & \geq \pr{\forall 0\leq t\leq T : X(t) \leq \phi(t)}\cdot (c \log T)^{-1/H}, \label{eqn:lowerboundwithout}
\end{align}
for sufficiently large $T$, since (using $\kappa>1$)
$$
\int_0^T e^{ \phi(u)} \, \d u = \int_0^{(\kappa \log T)^{1/H}} e^1 \, \d u + \int_{(\kappa \log T)^{1/H}}^T T^{-\kappa} \d u \leq 2 e (\kappa \log T)^{1/H}.
$$

By Slepian's inequality \cite{slepian1962} and noting that $\E X(t) X(s) \geq 0$, we have
\begin{align} \label{eqn:slepde}
& \pr{\forall 0\leq t\leq T : X(t) \leq \phi(t)}  \notag \\
 \geq~ & \pr{\sup_{0\leq t\leq (\kappa \log T)^{1/H}} X(t) \leq 1} \cdot \pr{\sup_{(\kappa \log T)^{1/H}\leq t\leq T} X(t)\leq -\kappa \log T}.
\end{align}
The first factor in (\ref{eqn:slepde}), by the lower bound in Theorem~\ref{thm:main} (which was proved in Section~\ref{sec:lower}), can be estimated as follows:
\begin{equation} \label{eqn:usemolchan}
\pr{\sup_{0\leq t\leq  (\kappa\log T)^{1/H}} X(t)\leq 1}  \geq \left[ (\kappa \log T)^{1/H}\right]^{-(1-H) + o(1)}.
\end{equation}
The second factor in (\ref{eqn:slepde}) equals
\begin{align*}
& \pr{\sup_{1 \leq s \leq T (\kappa \log T)^{-1/H}} X(s (\kappa \log T)^{1/H})\leq -\kappa \log T}\\
 =~& \pr{\sup_{1 \leq s \leq T  (\kappa \log T)^{-1/H}} (\kappa \log T) X(s)\leq -\kappa \log T}\\
 =~& \pr{\sup_{1 \leq s \leq T  (\kappa \log T)^{-1/H}} X(s)\leq -1}\\
 \geq ~& \pr{ \sup_{1 \leq s \leq T} X(s)\leq -1 }.
\end{align*}

Let $H\geq 1/2$. Then the increments of FBM are positively correlated. Therefore, by Slepian's lemma (second step),
\begin{align}
 & \pr{ \sup_{1 \leq s \leq T} X(s)\leq -1 } \notag \\
\geq ~ &\pr{ \sup_{1 \leq s \leq T} (X(s)-X(1)) \leq 1, X(1) \leq -2 } \notag \\
\geq ~ &\pr{ \sup_{1 \leq s \leq T} (X(s)-X(1)) \leq 1} \cdot \pr{ X(1) \leq -2 } \notag \\
= ~&\pr{ \sup_{1 \leq s \leq T} X(s-1) \leq 1} \cdot k. \notag \\
\geq ~&\pr{ \sup_{0 \leq s \leq T } X(s) \leq 1} \cdot k. \label{eqn:simpleesthlarger}
\end{align}
Putting these estimates together with Molchan's result (Statement~1 in~\cite{molchan1999}) for $I(T)$ yields the assertion:
\begin{align}
& c' T^{-(1-H)}  \geq I(T) \notag \\ & \geq \pr{ \sup_{0 \leq s \leq T} X(s) \leq 1} \cdot k \cdot (\log T)^{-(1-H)/H+o(1)} (c \log T)^{-1/H}. \tag*{$\square$}
\end{align}

In the case $H<1/2$, the proof is more involved, even though the main idea is the same. We start with the following purely analytic fact.

\begin{lem} \label{lem:analyticfact} 
Let $\ell(t):=2(\log \log (t e^{e}))^{\lambda}$ with $\lambda>0$. Let $0<\alpha\leq 1$. Then there is an $s_0=s_0(H)\geq 1$ such that, for all $t\geq s\geq s_0$,
\begin{align*}
 & \left( \ell(t)^2 t^{\alpha} - \ell(t) ( t^\alpha + 1 - (t-1)^\alpha) + 1 \right)^{1/\alpha}\\
  &~~~\geq~  \left( \ell(t)^2 t^{\alpha} - \ell(t)\ell(s) ( t^\alpha + s^\alpha - (t-s)^\alpha) + \ell(s)^2 s^\alpha \right)^{1/\alpha}  \\ &~~~~~~ ~~~ + \left( \ell(s)^2 s^{\alpha} - \ell(s) ( s^\alpha + 1 - (s-1)^\alpha) + 1 \right)^{1/\alpha}.
\end{align*}
\end{lem}

The proof of this elementary lemma is given in the appendix. We continue with an auxiliary lemma. In view of (\ref{eqn:simpleesthlarger}), this lemma highlights the difficulties with one-sided exit problems for general processes.

\begin{lem} \label{lem:hlesshalfdetour} Let $0<H<1/2$. Then there is an $s_0=s_0(H)\geq 1$ and a constant $k>0$ such that for any $K\in \R$
\begin{equation} \label{eqn:difficultsiesstimate}
\pr{ \sup_{ s_0\leq t \leq T} X(t) \leq -K} \geq \pr{ \sup_{0\leq t \leq k T (\log \log T)^{1/(4H)}}  X(t) \leq 1} (\log T)^{-o(1)}.
\end{equation}
\end{lem}

{\bf Proof.} 
Let $\ell(t):=2(\log \log (t e^{e}))^{\lambda}$ with $\lambda:=1/4$ and define
$$
Y(t):= \ell(t) X(t) -  X(1),\qquad t\geq 1.
$$
The idea of the proof is that $(Y(t))_{t\geq 1}$ and $X(1)$ are positively correlated (unlike $(X(t)-X(1))_{t\geq 1}$ and $X(1)$); but $Y(t)$ is essentially the same as $\ell(t) X(t)$, at least for large $t$.

Note that, for $t\geq 1$,
\begin{align}
\E Y(t) X(1) & =  \ell(t) \E X(t)  X(1) -  \E X(1)^2 \notag \\ & =  \ell(t)\frac{1}{2} (t^{2H} + 1 - (t-1)^{2H}) - 1 \geq  \ell(t)/2- 1 \geq 0. \label{eqn:fundonoR0}
\end{align}

Furthermore, define the function $f$ on $[1,\infty)$ by
\begin{align}
\E Y(t)^2 &= \ell(t)^2 \E X(t)^2 -2 \ell(t) \E X(t) X(1) + \E X(1)^2 \notag \\
         & = \ell(t)^2 t^{2H} -2 \ell(t) \frac{1}{2} (t^{2H}+1-(t-1)^{2H})  +  1 \notag \\
         & =: f(t)^{2H}. \label{eqn:defnf}
\end{align}
Then $f$ is increasing since $f(t)^{2H}$ is (as can be seen immediately by differentiating). In fact, for some constant $k>0$,
\begin{equation} \label{eqn:app}
f(t)\leq k (\log \log t)^{1/(4H)} t,\qquad\text{as $t\to\infty$.}
\end{equation}
Furthermore, the definition of $f$ in (\ref{eqn:defnf}) is such that
$$
\E Y(t)^2 = \E X(f(t))^2,\qquad t\geq 1.
$$
Further, one checks that for some $s_0=s_0(H)\geq 1$,
\begin{equation} \label{eqn:fundonoR1}
\E Y(t) Y(s) \geq \E X(f(t)) X(f(s)),\qquad t,s\geq s_0.
\end{equation}
Indeed, let $t\geq s\geq 1$ and recall that this is equivalent to 
\begin{equation} \label{eqn:fundono}
\E |Y(t) - Y(s)|^2 \leq \E |X(f(t)) - X(f(s))|^2.
\end{equation}
Note that (\ref{eqn:fundono}) can be rewritten as
$$
\ell(t)^2 t^{2H} - 2 \ell(t)\ell(s) \E X(t) X(s) + \ell(s)^2 s^{2H} \leq |f(t)-f(s)|^{2H}.
$$
This is Lemma~\ref{lem:analyticfact} with $\alpha=2H <1$.

Now we are ready for the main argument. We use Slepian's lemma together with (\ref{eqn:fundonoR0}) and  (\ref{eqn:fundonoR1}) in (\ref{eqn:ussl1}) and (\ref{eqn:ussl2}), respectively, to get that
\begin{align}
& \pr{ \sup_{ s_0\leq t \leq T} X(t) \leq -K} \notag \\
&= \pr{ \ell(t) X(t) \leq - \ell(t) K, \forall s_0\leq t \leq T} \notag \\
&\geq \pr{ \sup_{s_0\leq t \leq T} Y(t) \leq 1 , X(1) \leq - \ell(T)K -1}  \notag \\
&\geq  \pr{ \sup_{s_0\leq t \leq T}  Y(t) \leq 1} \cdot \pr{ X(1) \leq - \ell(T)K -1} \label{eqn:ussl1}\\
&\geq  \pr{ \sup_{s_0\leq t \leq T} X(f(t)) \leq 1 } \cdot \pr{ X(1) \leq - 4(\log \log (T e^e))^{1/4} K} \label{eqn:ussl2}.
\end{align}
The second term is of order $(\log T)^{-o(1)}$. The first term (by (\ref{eqn:app})) can be estimated from below by
\begin{align}
& \pr{ \sup_{ 0\leq t \leq T k (\log \log T)^{1/(4H)}} X(t) \leq 1}.  \tag*{$\square$}
\end{align}

Now we can prove the lower bound also in the case $H<1/2$.

\begin{proof}[ of the lower bound in Theorem~\ref{thm:main}, case $H<1/2$]
Let $s_0\geq 1$ be the constant from Lemma~\ref{lem:hlesshalfdetour}. As in (\ref{eqn:lowerboundwithout})-(\ref{eqn:usemolchan}), we obtain
$$
I(T) \geq  \pr{\sup_{(\kappa \log T)^{1/H}\leq t\leq T} X(t)\leq -\kappa \log T}\cdot (c \log T)^{-(2-H)/H-o(1)}.
$$
The first factor on the right-hand side equals:
\begin{align*}
& \pr{\sup_{s_0 \leq s \leq T s_0 (\kappa \log T)^{-1/H}} X(s s_0^{-1} (\kappa \log T)^{1/H})\leq -\kappa \log T}\\
 = ~& \pr{\sup_{s_0 \leq s \leq T  s_0 (\kappa \log T)^{-1/H}} s_0^{-1/H}  (\kappa \log T) X(s)\leq -\kappa \log T}\\
 =~& \pr{\sup_{s_0 \leq s \leq T s_0 (\kappa \log T)^{-1/H}} X(s)\leq -s_0^{1/H} }\\
 \geq ~& \pr{ \sup_{s_0 \leq s \leq T} X(s)\leq -s_0^{1/H} }.
\end{align*}
Using Lemma~\ref{lem:hlesshalfdetour}, this shows
$$
I(T) \geq \pr{ \sup_{ 0\leq t \leq T k (\log \log T)^{1/(4H)}} X(t) \leq 1} \cdot (c \log T)^{-(2-H)/H-o(1)}.
$$
Putting these estimates together with Molchan's result (Statement~1 in \cite{molchan1999}) for $I(T)$ yields the assertion.
\end{proof}


\begin{remark}
Let us comment on why it seems plausible that the lower bound could hold true without logarithmic loss factors.
If we choose
$$
\tilde\phi(t):=1-\kappa \log_+ t:=\begin{cases}
          1  & 0\leq t \leq 1,\\
	  1 - \kappa \log t & 1 \leq t \leq T
         \end{cases}
$$
instead of $\phi$ in the proof, we obtain as in (\ref{eqn:lowerboundwithout}):
\begin{equation}  \label{eqn:otherdrift}
c\, T^{-(1-H)} \geq I(T) \geq \pr{\sup_{0\leq t\leq T} (X(t) +\kappa \log_+ t) \leq 1}\cdot c',
\end{equation}
with some constants $c,c'>0$. In view of \cite{uchiyama1980}, it seems plausible that the latter probability has the same asymptotic rate (in the weak sense, but without additional logarithmic factors) as $\pr{\sup_{0\leq t\leq T} X(t) \leq 1}$. However, at the moment there seems to be no way of proving this.
\end{remark}

\begin{remark}
Note that starting from (\ref{eqn:otherdrift}) or (\ref{eqn:lowerboundwithout}), one immediately obtains Molchan's result (\ref{eqn:molchan}) by the application of Proposition~3.1 in \cite{aurzadadereich2009}. This presents a new and simple proof of the upper bound in (\ref{eqn:molchan}) for all $0<H<1$.

We remark that this argument works for any Gaussian process such that the function $g(t)=\kappa \log_+ t$ is bounded from above by some function from the RKHS of the process, cf.\ Proposition~3.1 in \cite{aurzadadereich2009}.
\end{remark}

\noindent {\bf Acknowledgements:} I would like to thank Michael Scheutzow (Berlin) for discussions on this subject. This work was supported by the DFG Emmy Noether programme.

\bibliographystyle{plain}

\appendix
\section*{Appendix}

\begin{proof}[ of Lemma~\ref{lem:analyticfact}]
 {\it Step 1:} We show that for all $t\geq s\geq 1$
$$
(\ell(t) (t-s)^\alpha - \ell(t) s^\alpha + \ell(s) s^\alpha )^{1/\alpha} \leq (\ell(t) t^\alpha - \ell(t) s^\alpha + \ell(s) s^\alpha )^{1/\alpha} - \ell(s)^{1/\alpha} s.
$$
To see this we rewrite the inequality as follows:
$$
(\ell(s) s)^{1/\alpha} \left(  (y z^\alpha - y + 1 )^{1/\alpha}- (y (z-1)^\alpha - y + 1 )^{1/\alpha} - 1 \right)\geq 0,
$$
where $y:=\ell(t)/\ell(s)\geq 1$, $z:=t/s\geq 1$, and $z\geq y$. In fact, one can show that
\begin{equation} \label{eqn:yetanotherhelp}
(y z^\alpha - y + 1 )^{1/\alpha}-(y (z-1)^\alpha - y + 1 )^{1/\alpha}\geq 1,\qquad \forall z\geq y\geq 1.
\end{equation}
Indeed, first one verifies that
\begin{equation} \label{eqn:indrad}
\frac{y(z^\alpha-1)+1}{y((z-1)^\alpha-1)+1} \geq \frac{z^\alpha}{(z-1)^\alpha},\qquad \forall z\geq y\geq 1,
\end{equation}
by observing that the left-hand side is increasing in $y$ (note that the function is of the type $y\mapsto (ya+1)/(yb+1)$ with $a\geq b>0$). Now note that the left-hand side of (\ref{eqn:yetanotherhelp}) equals
\begin{align*}
& (y (z^\alpha - 1) + 1 )^{1/\alpha}\left( 1 -\left(\frac{y ((z-1)^\alpha - 1) + 1}{y (z^\alpha - 1) + 1} \right)^{1/\alpha} \right)\\
\geq~ & (y (z^\alpha - 1) + 1 )^{1/\alpha}\left( 1 -\frac{z-1}{z} \right)\\
 =~& \left(\frac{y (z^\alpha - 1) + 1}{z^\alpha} \right)^{1/\alpha} \geq \left(\frac{(z^\alpha - 1) + 1}{z^\alpha} \right)^{1/\alpha}=1,
\end{align*}
where we used (\ref{eqn:indrad}) in the first step. This shows (\ref{eqn:yetanotherhelp}).

{\it Step 2:} We show that there is an $s_0=s_0(\alpha)$ such that for all $t\geq s\geq s_0$
$$
\ell(t)^{1/\alpha} (t^\alpha -1)^{1/\alpha} - \ell(s)^{1/\alpha}(s^\alpha-1)^{1/\alpha} \geq (\ell(t) t^\alpha - \ell(t) s^\alpha + \ell(s) s^\alpha )^{1/\alpha} - \ell(s)^{1/\alpha} s.
$$
To see this define the functions
$$
h_1(t):=\ell(t)^{1/\alpha} (t^\alpha -1)^{1/\alpha},\qquad h_2(t):=(\ell(t) t^\alpha - \ell(t) s^\alpha + \ell(s) s^\alpha )^{1/\alpha}.
$$
The assertion of Step 2 is that $h_1(t)-h_1(s)\geq h_2(t)-h_2(s)$. Since the functions $h_1, h_2$ are continuously differentiable, it is sufficient to show that $h_1'(t)\geq h_2'(t)$ for $t\geq s\geq s_0$.

We calculate:
$$
h_1'(t)=\frac{1}{\alpha} \ell(t)^{1/\alpha-1} \ell'(t)(t^\alpha -1)^{1/\alpha}+ \ell(t)^{1/\alpha} \frac{1}{\alpha}\, (t^\alpha -1)^{1/\alpha-1} \alpha t^{\alpha-1} 
$$
and
$$
h_2'(t)= \frac{1}{\alpha} (\ell(t) t^\alpha - \ell(t) s^\alpha + \ell(s) s^\alpha )^{1/\alpha-1}(\ell'(t) t^\alpha + \ell(t) \alpha t^{\alpha-1} - \ell'(t) s^\alpha).
$$
In order to see $h_1'(t)\geq h_2'(t)$ for $t\geq s\geq s_0$ we will show that
\begin{equation} \label{eqn:caseA}
 \ell(t)^{1/\alpha-1} \ell'(t)(t^\alpha -1)^{1/\alpha} \geq (\ell(t) t^\alpha - \ell(t) s^\alpha + \ell(s) s^\alpha )^{1/\alpha-1}(\ell'(t) t^\alpha  - \ell'(t) 2 )
\end{equation}
and
\begin{align} \label{eqn:caseB}
 & \ell(t)^{1/\alpha} (t^\alpha -1)^{1/\alpha-1} \alpha t^{\alpha-1} \notag  \\ \geq ~ & (\ell(t) t^\alpha - \ell(t) s^\alpha + \ell(s) s^\alpha )^{1/\alpha-1}(\ell(t) \alpha t^{\alpha-1} - \ell'(t) (s^\alpha-2)).
\end{align}

Let us first show (\ref{eqn:caseA}). First we note that it is nothing else but
$$
(t^\alpha -1)^{1/\alpha} \geq ( t^\alpha -  s^\alpha + \frac{\ell(s)}{\ell(t)} s^\alpha )^{1/\alpha-1}( t^\alpha  -  2 )
$$
Note that even
$$
(t^\alpha -1)^{1/\alpha} \geq ( t^\alpha+ 0 )^{1/\alpha-1}( t^\alpha  -  2 )
$$
holds for all $t$ sufficiently large. This finishes (\ref{eqn:caseA}).

Now we come to (\ref{eqn:caseB}). Dividing by $\ell(t)^{1/\alpha}$, (\ref{eqn:caseB}) reads as follows
\begin{equation} \label{eqn:vered}
(t^\alpha -1)^{1/\alpha-1} \alpha t^{\alpha-1} \geq ( t^\alpha -  s^\alpha (1- \frac{\ell(s)}{\ell(t)}) )^{1/\alpha-1}(\alpha t^{\alpha-1} - \frac{\ell'(t)}{\ell(t)} (s^\alpha-2)).
\end{equation}

{\it Step 2a:} If $(\log t )^2\leq (s^\alpha -2)/\alpha$ then we even have
$$
(t^\alpha -1)^{1/\alpha-1} \alpha t^{\alpha-1} \geq ( t^\alpha -  0 )^{1/\alpha-1}(\alpha t^{\alpha-1} - \frac{\ell'(t)}{\ell(t)} (s^\alpha-2)),
$$
i.e.\
$$
(( t^\alpha  )^{1/\alpha-1}- (t^\alpha -1)^{1/\alpha-1} )\alpha t^{\alpha-1} ( t^\alpha )^{-1/\alpha+1} \, \frac{\ell(t)}{\ell'(t)} \leq  (s^\alpha-2).
$$
This can be seen as follows. Note that the left-hand side of this inequality behaves asymptotically as $c (\log t) (\log \log t)$. So, for those $t$ with $(\log t)^2\leq (s^\alpha -2)/\alpha$ the inequality holds.

{\it Step 2b:} On the other hand, if  $(\log t )^2\geq (s^\alpha -2)/\alpha$ then $2 \log \log t \geq \log ((s^\alpha -2)/\alpha) \geq c_\alpha \log s$. Thus, $\ell(t)\geq c_\alpha' (\log s)^\lambda$ and therefore
$$
q:=s^\alpha (1- \frac{\ell(s)}{\ell(t)})\geq s^\alpha (1- \frac{\ell(s)}{c_\alpha' (\log s)^\lambda})  \geq 1,
$$
for sufficiently large $s$. This shows that the inequality (\ref{eqn:vered}) is also satisfied:
$$
(t^\alpha -1)^{1/\alpha-1} ( \alpha t^{\alpha-1} - 0) \geq ( t^\alpha -  q )^{1/\alpha-1}(\alpha t^{\alpha-1} - \frac{\ell'(t)}{\ell(t)} (s^\alpha-2))).
$$

{\it Step 3:} We show that there is an $s_0=s_0(\alpha)$ such that for all $t\geq s\geq s_0$
$$
\ell(t)^{1/\alpha} (t^\alpha -1)^{1/\alpha} - \ell(s)^{1/\alpha}(s^\alpha-1)^{1/\alpha} - (\ell(t) (t-s)^\alpha - \ell(t) s^\alpha + \ell(s) s^\alpha )^{1/\alpha} \geq 0.
$$
This follows directly from Steps~1 and~2:
\begin{align*}
& \ell(t)^{1/\alpha} (t^\alpha -1)^{1/\alpha} - \ell(s)^{1/\alpha}(s^\alpha-1)^{1/\alpha} \\
\geq ~& (\ell(t) t^\alpha - \ell(t) s^\alpha + \ell(s) s^\alpha )^{1/\alpha} - \ell(s)^{1/\alpha} s\\
\geq ~& (\ell(t) (t-s)^\alpha - \ell(t) s^\alpha + \ell(s) s^\alpha )^{1/\alpha}.
\end{align*}

{\it Step 4:} We finally show the assertion. In the following calculation we use in the first, third, and fourth step that
$$
(x+z)^{1/\alpha}-(y+z)^{1/\alpha} \geq x^{1/\alpha}- y^{1/\alpha},\qquad \forall x\geq y\geq 0, z\geq 0.
$$
Then we obtain:
\begin{align*}
 & \left( \ell(t)^2 t^{\alpha} - \ell(t) ( t^\alpha + 1 - (t-1)^\alpha) + 1 \right)^{1/\alpha}\\
 & - \left( \ell(s)^2 s^{\alpha} - \ell(s) ( s^\alpha + 1 - (s-1)^\alpha) + 1 \right)^{1/\alpha} \\ 
 & - \left( \ell(t)^2 t^{\alpha} - \ell(t)\ell(s) ( t^\alpha + s^\alpha - (t-s)^\alpha) + \ell(s)^2 s^\alpha \right)^{1/\alpha} \\
\geq ~& \left( \ell(t)^2 t^{\alpha} - \ell(t)  \right)^{1/\alpha} \\
 & - \left( \ell(s)^2 s^{\alpha} - \ell(s) ( s^\alpha + 1 - (s-1)^\alpha) +\ell(t)(t^\alpha -(t-1)^\alpha)\right)^{1/\alpha} \\ 
 & - \left( \ell(t)^2 t^{\alpha} - \ell(t)\ell(s) ( t^\alpha + s^\alpha - (t-s)^\alpha) + \ell(s)^2 s^\alpha \right)^{1/\alpha}\\
\geq ~& \left( \ell(t)^2 t^{\alpha} - \ell(t)  \right)^{1/\alpha} \\
 & - \left( \ell(s)^2 s^{\alpha} - \ell(s)\right)^{1/\alpha} \\ 
 & - \left( \ell(t)^2 t^{\alpha} - \ell(t)\ell(s) ( t^\alpha + s^\alpha - (t-s)^\alpha) + \ell(s)^2 s^\alpha \right)^{1/\alpha}\\
\geq ~& \left( \ell(t) \ell(s) t^{\alpha} - \ell(t)  \right)^{1/\alpha} \\
 & - \left( \ell(s)^2 s^{\alpha} - \ell(s)\right)^{1/\alpha} \\ 
 & - \left(  - \ell(t)\ell(s) s^\alpha + \ell(t)\ell(s) (t-s)^\alpha + \ell(s)^2 s^\alpha \right)^{1/\alpha}\\
\geq ~& \left( \ell(t) \ell(s) t^{\alpha} - \ell(t) \ell(s)  \right)^{1/\alpha} \\
 & - \left( \ell(s)^2 s^{\alpha} - \ell(s)+\ell(t) - \ell(t)\ell(s)\right)^{1/\alpha} \\ 
 & - \left(  - \ell(t)\ell(s) s^\alpha + \ell(t)\ell(s) (t-s)^\alpha + \ell(s)^2 s^\alpha \right)^{1/\alpha}\\
\geq ~& \left( \ell(t) \ell(s) t^{\alpha} - \ell(t) \ell(s)  \right)^{1/\alpha} \\
 & - \left( \ell(s)^2 s^{\alpha} - \ell(s)^2\right)^{1/\alpha} \\ 
 & - \left(  - \ell(t)\ell(s) s^\alpha + \ell(t)\ell(s) (t-s)^\alpha + \ell(s)^2 s^\alpha \right)^{1/\alpha},
\end{align*}
since $z:=-\ell(t)(t^\alpha -(t-1)^\alpha)+1\geq 0$ (first step), $\ell(s) ( s^\alpha - (s-1)^\alpha) -\ell(t)(t^\alpha -(t-1)^\alpha)\geq 0$ (second step), $z:=\ell(t)^2 t^{\alpha} - \ell(t)\ell(s) t^\alpha\geq 0$ (third step), $z:=\ell(t)\ell(s)-\ell(t)\geq 0$ (fourth step), and $- \ell(s)+\ell(t)-\ell(t)\ell(s)\leq -\ell(s)^2$ (fifth step).

Now, the term can be divided by $\ell(s)^{1/\alpha}$; and it thus remains to be seen that
$$
\ell(t)^{1/\alpha} (  t^{\alpha} - 1  )^{1/\alpha} - \ell(s)^{1/\alpha} (  s^{\alpha} -1)^{1/\alpha}  - ( \ell(t) (t-s)^\alpha - \ell(t) s^\alpha + \ell(s) s^\alpha )^{1/\alpha} \geq 0.
$$
However, this is exactly what was asserted in Step 3.
\end{proof}

\end{document}